\documentclass[oneside]{amsart}
\newcommand{\formatswitch}{paper}

\usepackage{amssymb}
\usepackage{ifthen}
\usepackage{verbatim}
\usepackage{epsfig}
\usepackage[all]{xypic}

\ifthenelse{\equal{\formatswitch}{big}}{
\setlength{\textwidth}{432pt}
\setlength{\oddsidemargin}{18.8775pt}

\setcounter{tocdepth}{1}
}{
\ifthenelse{\equal{\formatswitch}{paper}}{

\setcounter{tocdepth}{1}
}{
\setcounter{tocdepth}{2}
}
}
\DeclareMathAlphabet\EuScript{U}{eus}{m}{n}
\DeclareMathAlphabet\EuScriptb{U}{eus}{b}{n}

\newcommand{\claimenum}{\renewcommand{\theenumi}{\alph{enumi}}
 \renewcommand{\labelenumi}{\textit{(\theenumi)}}
 \renewcommand{\theenumii}{\roman{enumii}}
 \renewcommand{\labelenumii}{\textit{(\theenumii)}}
 \begin{enumerate}}
\newcommand{\claimenumend}{\end{enumerate}}

\newcommand{\romanenum}{\renewcommand{\theenumi}{\roman{enumi}}
 \renewcommand{\labelenumi}{\textit{(\theenumi)}}
 \renewcommand{\theenumii}{\alph{enumii}}
 \renewcommand{\labelenumii}{\textit{(\theenumii)}}
 \begin{enumerate}}
\newcommand{\romanenumend}{\end{enumerate}}

\newtheorem{dummy}{realdumb}
\newtheorem{thm}{Theorem}

\newtheorem{lemma}[dummy]{Lemma}

{\theoremstyle{definition} }

{\theoremstyle{definition} }

\renewcommand{\text}{\mathrm}

\newcommand{\strutdepth}{\dp\strutbox}
\newcommand{\marginalnote}[1]
   {\strut\vadjust{\kern-\strutdepth\domarginalnote{#1}}}
\newcommand{\domarginalnote}[1]{\vtop to \strutdepth{
  \baselineskip\strutdepth
   \vss\llap{ #1\ \ }\null}}  
\newcounter{showlabelflag}
\newcounter{makelabelflag}
\newcommand{\showlabels}{\setcounter{showlabelflag}{1}}
\newcommand{\nolabels}{\setcounter{makelabelflag}{2}}
\newcommand{\makelabels}{\setcounter{makelabelflag}{1}}

\newcommand{\mylabel}[1]{
  \ifthenelse{\value{makelabelflag}=1}
    {\label{#1}}{}
  \ifthenelse{\value{showlabelflag}=1}
    {\marginpar{#1}}{}\relax}

\ifthenelse{\equal{\formatswitch}{draft}}{
\showlabels
}{\relax}

\begin{document}

\bibliographystyle{amsplain}

\nolabels
\begin{center}{\bfseries On the Baker's map and  
the Simplicity \\ of the Higher Dimensional Thompson Groups \(nV\)
}\end{center}
\vspace{3pt}
\begin{center}{MATTHEW G. BRIN}\end{center}
\vspace{4pt}

\noindent {\bfseries Abstract}.  {\itshape We show that the baker's
map is a product of transpositions (particularly pleasant
involutions), and conclude from this that an existing very short
proof of the simplicity of Thompson's group \(V\) applies with equal
brevity to the higher dimensional Thompson groups \(nV\).}

\vspace{4pt}

\ifthenelse{\equal{\formatswitch}{paper}}{
\relax
}{
\tableofcontents
\clearpage
}

\makelabels
\ifthenelse{\equal{\formatswitch}{draft}}{
\showlabels
}{\relax}


Of the original groups \(F\subseteq T\subseteq V\) of Richard
Thompson (see \cite{CFP}), all are infinite and finitely presented,
and the last two are simple.  An infinite family of groups \(nV\),
\(n\in \{1, 2, 3, \ldots, \omega\}\), of which \(1V=V\), is
introduced in \cite{brin:hd3} where it is shown that \(2V\) is
infinite and simple and not isomorphic to \(V\).  A finite
presentation for \(2V\) is given in \cite{brin:hd4}, it is shown
that \(nV\) and \(mV\) are isomorphic only when \(m=n\) in
\cite{bleak+lanoue}, and metric properties of \(2V\) are studied in
\cite{burillo:hdmetric}.

A very short argument that \(V=1V\) is generated by transpositions is
given in Section 12 of \cite{brin:hd3}, followed by an equally short
argument based on this fact (due to Rubin) that \(V\) is simple.  It
is also shown in that section that the baker's map (an element of
\(2V\)) prevents the first argument from showing that \(2V\) is
generated by transpositions.  As a result, the proof in \cite{brin:hd3}
of the simplicity of \(2V\) is rather involved and is based on
calculations which show that the abelianization of \(2V\) is
trivial.

Here we give a short proof that the baker's maps in \(2V\) are
products of transpositions in \(2V\).  From there it is an easy
exercise to combine this with the material in Section 12 of
\cite{brin:hd3} to give a short proof of the simplicity of \(2V\)
and also to extend all the results to all of the \(nV\), \(n\le
\omega\).

Longer arguments for simplicity exist. Presentations (finite when
\(n<\omega\)) for the \(nV\), \(n\le\omega\), are given in
\cite{matucci+hennig}, and one can calculate from these
presentations that each of the \(nV\), \(n \le \omega\), has trivial
abelianization.  From the arguments in Section 3 of \cite{brin:hd3}
(which are about \(2V\) but, as noted in 4.1 of \cite{brin:hd3},
they apply as well to the \(nV\), \(n \le \omega\)) it then follows
that each of the \(nV\), \(n \le \omega\), is simple.

To keep this paper brief, we use notation, terminology and graphics
from \cite{brin:hd3}, and from this point we assume that the reader
is familiar with their meanings.

\begin{lemma} Any baker's map in \(2V\) is a product of finitely
many proper transpositions from \(2V\).  \end{lemma}

The (primary) baker's map is given by the following.

\[
\xy
(-9,-9); (-9,9)**@{-}; (9,9)**@{-}; (9,-9)**@{-}; (-9,-9)**@{-};
(0,-9); (0,9)**@{-};
(-4.5,0)*{\scriptstyle0};
(4.5,0)*{\scriptstyle1};
\endxy
\quad
\longrightarrow
\quad
\xy
(-9,-9); (-9,9)**@{-}; (9,9)**@{-}; (9,-9)**@{-}; (-9,-9)**@{-};
(-9,0); (9,0)**@{-};
(0,-4.5)*{\scriptstyle0};
(0,4.5)*{\scriptstyle1};
\endxy
\]

A (secondary) baker's map is given by a pair of patterns that are
identical and identically numbered with one exception: for one
singly divided rectangle in the domain and for the corresponding
singly divided rectangle in the range, the division is vertical in
the domain and horizontal in the range.  An example is below.

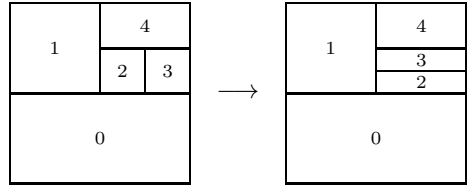
\[
\xy (0,0); (0,24)**@{-}; (24,24)**@{-}; (24,0)**@{-};
(0,0)**@{-};
(0,12); (24,12)**@{-}; (12,12); (12,24)**@{-};
(12,18); (24,18)**@{-};
(18,12); (18,18)**@{-};
(15,15)*{\scriptstyle2};
(21,15)*{\scriptstyle3};
(18,21)*{\scriptstyle4};
(6,18)*{\scriptstyle1};
(12,6)*{\scriptstyle0};
\endxy
\quad
\xy
(0,0)*{\ };
(0,12)*{\longrightarrow};
\endxy
\quad
\xy (0,0); (0,24)**@{-}; (24,24)**@{-}; (24,0)**@{-};
(0,0)**@{-};
(0,12); (24,12)**@{-}; (12,12); (12,24)**@{-};
(12,18); (24,18)**@{-};
(12,15); (24,15)**@{-};
(18,13.5)*{\scriptstyle2};
(18,16.5)*{\scriptstyle3};
(18,21)*{\scriptstyle4};
(6,18)*{\scriptstyle1};
(12,6)*{\scriptstyle0};
\endxy
\]

We refer to the rectangle containing the non-identity part of the
baker's map as the {\itshape support} of the baker's map.

A transposition is given by a pair of identical
patterns that are numbered identically except for a switch of two of
the numbers.  The transposition is {\itshape proper} if there are
more than two rectangles in each pattern.

In the following, we say that two elements are identical modulo
transpositions if each is a product of elements and the two
products are identical if the transpositions are removed.

1. {\itshape Any baker's map \(\beta\) is, modulo transpositions, a
product of a baker's map whose support is the left half of the
support of \(\beta\), with a baker's map whose support is the right
half of the support of \(\beta\).}  In the pictures below, we show
only the support of \(\beta\).  The first arrow is the baker's map
\(\beta\) in a reducible form.  The second is a proper
transposition.  The composition is the promised product of two
baker's maps.

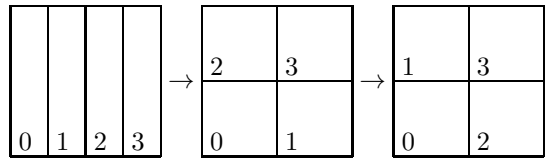
\[
\xy
(-10,-10); (-10,10)**@{-}; (10,10)**@{-}; (10,-10)**@{-}; (-10,-10)**@{-};
(0,-10); (0,10)**@{-}; (-8,-8)*{0}; (2,-8)*{2};
(-5,-10); (-5,10)**@{-}; (-3,-8)*{1}; (7,-8)*{3};
(5,-10); (5,10)**@{-}; 
\endxy
\rightarrow
\xy
(-10,-10); (-10,10)**@{-}; (10,10)**@{-}; (10,-10)**@{-}; (-10,-10)**@{-};
(-10,0); (10,0)**@{-}; (-8,-8)*{0}; (-8,2)*{2};
(0,-10); (0,10)**@{-}; (2,-8)*{1};  (2,2)*{3};
\endxy
\rightarrow
\xy
(-10,-10); (-10,10)**@{-}; (10,10)**@{-}; (10,-10)**@{-}; (-10,-10)**@{-};
(-10,0); (10,0)**@{-}; (-8,-8)*{0}; (-8,2)*{1};
(0,-10); (0,10)**@{-}; (2,-8)*{2};  (2,2)*{3};
\endxy
\]

2. {\itshape Any baker's map \(\beta\) is, modulo transpositions, a
product of a baker's map whose support is the upper half of the
support of \(\beta\), with a baker's map whose support is the lower
half of the support of \(\beta\).}  The relevant pictures follow and
the comments are as in 1.

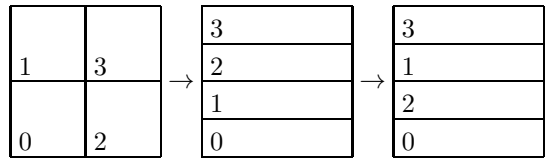
\[
\xy
(-10,-10); (-10,10)**@{-}; (10,10)**@{-}; (10,-10)**@{-}; (-10,-10)**@{-};
(-10,0); (10,0)**@{-}; (-8,-8)*{0}; (-8,2)*{1};
(0,-10); (0,10)**@{-}; (2,-8)*{2};  (2,2)*{3};
\endxy
\rightarrow
\xy
(-10,-10); (-10,10)**@{-}; (10,10)**@{-}; (10,-10)**@{-}; (-10,-10)**@{-};
(-10,0); (10,0)**@{-}; (-8,-8)*{0}; (-8,-3)*{1};
(-10,-5); (10,-5)**@{-}; (-8,2)*{2};  (-8,7)*{3};
(-10,5); (10,5)**@{-}; 
\endxy
\rightarrow
\xy
(-10,-10); (-10,10)**@{-}; (10,10)**@{-}; (10,-10)**@{-}; (-10,-10)**@{-};
(-10,0); (10,0)**@{-}; (-8,-8)*{0}; (-8,-3)*{2};
(-10,-5); (10,-5)**@{-}; (-8,2)*{1};  (-8,7)*{3};
(-10,5); (10,5)**@{-}; 
\endxy
\]

In the following, ``arbitrarily small'' means having support with
diameter smaller than an arbitrarily chosen positive real.

3. {\itshape Any baker's map is, modulo transpositions, a product of
arbitrarily small baker's maps.}  This follows from 1 and 2.

4. {\itshape A product of a baker's map and an inverse of a baker's
map with disjoint supports is a product of transpositions.}  Let \(A\)
and \(B\) be the disjoint supports.  We refer to the rectangles in
figures below to describe a sequence of transpositions.  (a) Switch
\(A_0\) with \(B_0\).  (b) Switch \(A_1\) with \(B_1\).  (c) Switch
\(A\) with \(B\).  The composition of (a) with (b) with (c) in that
order is the desired result.

\[
A=\xy
(-9,-9); (-9,9)**@{-}; (9,9)**@{-}; (9,-9)**@{-}; (-9,-9)**@{-};
(0,-9); (0,9)**@{-};
(-4.5,0)*{A_0};
(4.5,0)*{A_1};
\endxy
\quad
,
\quad
B=\xy
(-9,-9); (-9,9)**@{-}; (9,9)**@{-}; (9,-9)**@{-}; (-9,-9)**@{-};
(-9,0); (9,0)**@{-};
(0,-4.5)*{B_0};
(0,4.5)*{B_1};
\endxy
\]

5. {\itshape If \(R\) is a rectangle in a pattern so that neither
  side of \(R\) has length more than \(\frac12\), then the baker's
  map with support \(R\) is a product of transpositions.}  The
assumptions make \(R\) one half of a rectangle \(A\) that is not all
of the unit square.  Thus there is a rectangle \(B\) that is
disjoint from \(A\).  Let \(S\) be the rectangle that is the ``other
half'' of \(A\).  Let \(\alpha\) be a product of a baker's map on
\(A\) with an inverse of a baker's map on \(B\).  By 4, this is a
product of transpositions.  By 1 or 2 we can modify \(\alpha\) so that
it is still a product of transpositions, and is a baker's map on each
of \(R\) and \(S\) and an inverse of a baker's map on \(B\).  Let
\(\beta\) be a product of a baker's map on \(B\) and an inverse of a
baker's map on \(S\).  By 4 this is a product of transpositions.  Now
the composition of \(\alpha\) with \(\beta\) gives the desired
result.

The lemma follows from 3 and 5. As discussed, this implies the
following.

\begin{thm} The \(nV\), \(n\le \omega\), are generated by
transpositions and are simple.  \end{thm}


\begin{thebibliography}{1}

\bibitem{bleak+lanoue}
Collin Bleak and Daniel Lanoue, \emph{A family of non-isomorphism results},
  ArXiv preprint: http://front.math.ucdavis.edu/0807.4955, 2008.

\bibitem{brin:hd3}
Matthew~G. Brin, \emph{Higher dimensional {T}hompson groups}, Geom. Dedicata
  \textbf{108} (2004), 163--192. \MR{MR2112673}

\bibitem{brin:hd4}
\bysame, \emph{Presentations of higher dimensional {T}hompson groups}, J.
  Algebra \textbf{284} (2005), no.~2, 520--558. \MR{MR2114568 (2007e:20062)}

\bibitem{burillo:hdmetric}
Jose Burillo, \emph{Metric properties of higher-dimensional {T}hompson's
  groups}, ArXiv preprint: http://front.math.ucdavis.edu/0810.3926, 2008.

\bibitem{CFP}
J.~W. Cannon, W.~J. Floyd, and W.~R. Parry, \emph{Introductory notes on
  {R}ichard {T}hompson's groups}, Enseign. Math. (2) \textbf{42} (1996),
  no.~3-4, 215--256. \MR{98g:20058}

\bibitem{matucci+hennig}
Johanna Hennig and Francesco Matucci, \emph{Presentations for the higher
  dimensional {T}hompson's groups}, preprint, 2009.

\end{thebibliography}

\providecommand{\bysame}{\leavevmode\hbox to3em{\hrulefill}\thinspace}
\providecommand{\MR}{\relax\ifhmode\unskip\space\fi MR }
\providecommand{\MRhref}[2]{%
  \href{http://www.ams.org/mathscinet-getitem?mr=#1}{#2}
}
\providecommand{\href}[2]{#2}

\noindent Binghamton University, Binghamton, NY 13902-6000, USA

\end{document}